\documentclass{amsart}
\usepackage[all]{xy}
\usepackage{amssymb}
\newcommand{\cC}{{\mathcal C}}
\newcommand{\cA}{{\mathcal A}}
\newcommand{\fA}{{\mathfrak A}}
\newcommand{\cD}{{\mathcal D}}
\newcommand{\cM}{{\mathcal M}}
\newcommand{\be}{{\bf 1}}

\newcommand{\BZ}{{\mathbb Z}}
\newcommand{\BN}{{\mathbb N}}
\newcommand{\BC}{{\mathbb C}}
\newcommand{\eps}{{\varepsilon}}
\newcommand{\sq}{$\square$}
\newcommand{\bi}{\bar \imath}
\newcommand{\bj}{\bar \jmath}
\newcommand{\Ve}{\mbox{Vec}}
\newcommand{\Rep}{\mbox{Rep}}
\newcommand{\Fun}{\mbox{Fun}}
\newcommand{\Hom}{\mbox{Hom}}
\newcommand{\Mod}{\mbox{Mod}}
\newcommand{\Bimod}{\mbox{Bimod}}
\newcommand{\Ind}{\mbox{Ind}}
\newcommand{\id}{\mbox{id}}
\newcommand{\End}{\mbox{End}}
\newcommand{\iHom}{\underline{\mbox{Hom}}}
\title{Module categories, weak Hopf algebras and modular invariants}
\author{Viktor Ostrik}
\email{ostrik@math.mit.edu}
\address{Department of Mathematics, MIT, 77 Massachusetts Ave., Cambridge,
MA 02139}
\thanks{The author was partially supported by NSF grant DMS-0098830}
\date{November 2001}
\begin{document}
\begin{abstract}We develop abstract nonsense for module categories over
monoidal categories (this is a straightforward categorification of 
modules over rings). As applications we show that any semisimple monoidal
category with finitely many simple objects is equivalent to the category of
representations of a weak Hopf algebra (theorem of T.~Hayashi) and classify 
module categories over the fusion category of $\widehat{sl}(2)$ at a 
positive integer level where we meet once again ADE classification pattern.
\end{abstract} 
\maketitle
\section{Introduction}
One considers the notion of an (abelian) monoidal category as a 
categorification
of the notion of a ring. From this point of view it is natural to define a 
{\em module category} over a monoidal category as a categorification of the 
notion of the module over a ring. Such a definition was given by I.~Bernstein, 
and L.~Crane and I.~B.~Frenkel, see \cite{B,CF}. The main point of this paper 
is to show that module categories is an extremely convenient language. 
Moreover, this notion is implicitly present in recent developments of
such subjects as (1) Boundary Conformal Field Theory,
see \cite{BPPZ, FS}; (2) Subfactors Theory, see \cite{BEK, O}; (3) Theory
of weak Hopf algebras, see \cite{NV, Sl}; (4) Theory of extensions of
vertex algebras, see \cite{KO}. My own motivation to study this notion comes
from the theory of affine Hecke algebras, see \cite{BO}.

The aim of this paper is to give basic definitions (sections 2 and 3),
to present some explanations of the relations with the subjects above (sections
4 and 5) and to give some examples (section 6). Perhaps this paper does not
contain new results, but I hope that its point of view, language and some
proofs are new.

This paper owes much to many people. I am grateful to Alexander Kirillov, Jr. 
who explained to me how to do calculations in tensor categories, and for 
collaboration in paper \cite{KO} which strongly influenced this paper; also 
Lemma 7 below is due to him. I am greatly indebted to Dmitri Nikshych who
taught me everything I know about weak Hopf algebras and patiently answered
my questions; Lemma 6 below is due to him. My sincere gratitude is due
to Pavel Etingof and David Kazhdan for many discussions that significantly 
clarified my understanding of the subject. Thanks are also due to Leonid 
Vainerman for bringing reference \cite{H} to my attention.

\section{Module categories}

\subsection{Based rings and modules}
We will follow the terminology borrowed from \cite{EK} and \cite{Lu}.
In what follows all rings, algebras are assumed to be associative with unit.
Let $\BZ_+$ be the set of nonnegative integers.

{\bf Definition 1.} (i) A {\em $\BZ_+-$basis} of an algebra free as a module
over $\BZ$ is
a basis $B=\{ b_i\}$ such that $b_ib_j=\sum_kc_{ij}^kb_k,\; c_{ij}^k\in \BZ_+$.

(ii) A {\em $\BZ_+-$ring} is an algebra over $\BZ$ with unit endowed with a
fixed $\BZ_+-$basis.

(iii) A {\em $\BZ_+-$module} over a $\BZ_+-$ring $A$ is $\BZ-$free 
$A-$module $M$ endowed with a fixed basis $\{ m_j\}$ such that
$b_im_j=\sum_kd_{ij}^km_k,\; d_{ij}^k\in \BZ_+$.

(iv) Two $\BZ_+-$modules $M_1, M_2$ over $A$ with bases $\{ m_i^1\}_{i\in I},\;
\{ m_j^2\}_{j\in J}$ are {\em equivalent} if and only if there exists a
bijection $\phi : I\to J$ such that the induced $\BZ-$linear map $\tilde
\phi$ of abelian groups $M_1, M_2$ defined by $\tilde \phi (m_i^1)=
m_{\phi (j)}^2$ is an isomorphism of $A-$modules.
 
(v) The {\em direct sum} of two $\BZ_+-$modules $M_1, M_2$ over $A$ is the
module $M_1\oplus M_2$ over $A$ with the basis being the union of the bases
of $M_1$ and $M_2$.

(vi) A $\BZ_+-$module $M$ over $A$ is {\em indecomposable} if it is not
equivalent to a direct of two nonzero $\BZ_+-$modules.

(vii) A {\em $\BZ_+-$submodule} of a $\BZ_+-$module $M$ over $A$ with basis
$\{ m_i\}_{i\in I}$ is an subset $J\subset I$ such  that abelian subgroup
of $M$ generated by $\{ m_i\}_{i\in J}$ is $A-$submodule.

(viii) A $\BZ_+-$module $M$ over $A$ is {\em irreducible} if any 
$\BZ_+-$submodule of $M$ is 0 or $M$.

{\bf Proposition 1.} (cf. \cite{Ga}) For a given $\BZ_+-$ring $A$ of finite 
rank over $\BZ$ there exist only finitely many irreducible inequivalent 
$\BZ_+-$modules over $A$.

{\bf Proof.} First of all it is clear that an irreducible $\BZ_+-$module $M$ 
over $A$ is of finite rank over $\BZ$. Let $\{ m_i\}_{i\in I}$ be the basis of 
$M$. Let us consider an element $b:=\sum_{b_i\in B}b_i$ of $A$. Let $b^2=
\sum_in_ib_i$ and let $N=\max_{b_i\in B}n_i$ ($N$ exists since $B$ is finite).
For any $i\in I$ let $bm_i=\sum_{k\in I}d_i^km_k$ and let $d_i=\sum_{k\in I}
d_i^k$. Let $i_0\in I$ be such that $d=d_{i_0}=\min_{i\in I}d_i$. Let
$b^2m_{i_0}=\sum_{i\in I}c_im_i$. Calculating $b^2m_{i_0}$ in two ways ---
as $(b^2)m_{i_0}$ and as $b(bm_{i_0})$ we have:
$$Nd\ge \sum_ic_i\ge d^2$$
and consequently $d\le N$. So there are only finitely many possibilities for
$|I|$, values of $c_i$ and consequently for expansions $b_im_k$ (since each
$m_k$ appears in $bm_{i_0}$). The Proposition is proved. \sq  

{\bf Definition 2.} (i) A $\BZ_+-$ring $A$ with basis $\{ b_i\}_{i\in I}$ is
called {\em based ring} if the following conditions hold

(a) There exists subset $I_0\subset I$ such that $1=\sum_{i\in I_0}b_i$.

(b) Let $\tau : A\to \BZ$ be the group homomorphism defined by
$$
\tau (b_i)=\left\{ \begin{array}{ccc}1&\mbox{if}& i\in I_0\\ 0&\mbox{if}& 
i\not \in I_0\\ \end{array}\right.
$$
There exists an involution $i\mapsto \bi$ of $I$ such that induced map
$a=\sum_{i\in I}a_ib_i\mapsto \bar a=\sum_{i\in I}a_ib_{\bi},\; a_i\in \BZ$
is an anti-involution of ring $A$ and such that 
$$
\tau (b_ib_j)=\left\{ \begin{array}{ccl}1&\mbox{if}& i=\bj\\ 0&\mbox{if}& 
i\ne \bj.\\ \end{array}\right.
$$

(ii) A {\em based module} over a based ring $A$ with basis $\{ b_i\}_{i\in I}$
is a $\BZ_+-$module $M$ with basis $\{ m_j\}_{j\in J}$ over $A$ such that 
$d_{ij}^k=d_{\bi k}^j$ where numbers $d_{ij}^k$ are defined in Definition 1
(iii).

(iii) A {\em unital based ring} is a based ring $A$ such that the set $I_0$
consists of one element.

{\bf Remark 1.} (i) It follows easily from definition that $i, j\in I_0,\; 
i\ne j$ implies that $b_i^2=b_i,\; b_ib_j=0,\; \bi=i$.

(ii) It is easy to see that for a given $\BZ_+-$ring $A$ being based ring is
a {\em property}, not additional structure.

(iii) A different terminology is used in physical literature: our notion
of a unital based ring corresponds to the notion of a {\em fusion rules
algebra} (at least in the commutative case) and the notion of a based module
corresponds to the notion of a NIM-rep, see \cite{BPPZ, FS}.

{\bf Lemma 1.} Let $M$ be a based module over a based ring $A$. If $M$ is
indecomposable as $\BZ_+-$module over $A$ then $M$ is irreducible as a
$\BZ_+-$module over $A$.

{\bf Proof.} Let $\{ m_i\}_{i\in I}$ be the basis of $M$.
By the definition of a based module the scalar product on $M$ 
defined by $(m_i,m_j)=\delta_{ij}$ is invariant with respect to the 
antiinvolution $a\mapsto \bar a$. Hence the orthogonal complement to an
$A-$submodule is again an $A-$submodule. Finally an orthogonal complement to 
a $\BZ-$submodule generated by $m_k, \; k\in K\subset I$ is the $\BZ-$submodule
generated by $m_j,\; j\in I-K$ and therefore the orthogonal complement to a 
based $A-$submodule is again a based $A-$submodule. \sq

\subsection{Monoidal categories} In this paper
we will consider only abelian semisimple categories over a field $k$
with finite dimensional Hom-spaces. If otherwise is not stated explicitly
we will assume that the field $k$ is algebraically closed.
All functors are assumed to be additive. 

{\bf Definition 3.} (see e.g. \cite{BK}) A {\em monoidal category} consists of
the following data: category $\cC$, functor $\otimes : \cC \times \cC \to \cC$,
functorial isomorphisms $a_{X,Y,Z}: (X\otimes Y)\otimes Z\to X\otimes
(Y\otimes Z)$, unit object $\be \in \cC$, functorial isomorphisms $r_X: 
X\otimes \be \to X$ and $l_X: \be \otimes X\to X$ subject to the following 
axioms:

1) Pentagon axiom: the diagram
$$
\xymatrix{&((X\otimes Y)\otimes Z)\otimes T \ar[dl]_{a_{X,Y,Z}\otimes id} 
\ar[dr]^{a_{X\otimes Y,Z,T}}&\\
(X\otimes (Y\otimes Z))\otimes T \ar[d]^{a_{X,Y\otimes Z,T}}&&(X\otimes Y)
\otimes (Z\otimes T) \ar[d]_{a_{X,Y,Z\otimes T}}\\X\otimes ((Y\otimes Z)
\otimes T) \ar[rr]^{id\otimes a_{Y,Z,T}}&&X\otimes (Y\otimes (Z\otimes T))}
$$
commutes.

2) Triangle axioms: the diagram
$$
\xymatrix{(X\otimes \be )\otimes Y\ar[rr]^{a_{X,\be ,Y}} \ar[dr]^{r_X\otimes 
id}&&X\otimes (\be \otimes Y)
\ar[dl]_{id\otimes l_Y}\\ &X\otimes Y&}
$$
commutes.

In what follows we will omit from notations associativity and unit isomorphisms
what is justified by Maclane coherence theorem, see e.g. \cite{BK}.
 
The Grothendieck group $K_0(\cC )$ of a monoidal category $\cC$ is endowed 
with the structure of $\BZ_+-$ring: multiplication is induced by $\otimes$, 
and $\BZ_+-$basis consists of classes of simple objects. We will say that 
a monoidal category $\cC$ is a {\em categorification} of $\BZ_+-$algebra 
$K_0(\cC )$. There exist examples when given $\BZ_+-$ring admits non unique 
categorification and when $\BZ_+-$ring admits no categorifications, see e.g. 
\cite{TY}.

{\bf Definition 4.} (see \cite{BK}) (i) Let $\cC$ be a monoidal category
and $X$ be an object in $\cC$. A {\em right dual} to $X$ is an object $X^*$
with two morphisms
$$
e_X: X^*\otimes X\to \be,\; i_X: \be \to X\otimes X^*
$$
such that the compositions
$$ X\stackrel{i_X\otimes id}{\longrightarrow}X\otimes X^*\otimes X
\stackrel{id \otimes e_X}{\longrightarrow}X$$
$$ X^*\stackrel{id\otimes i_X}{\longrightarrow}X^*\otimes X\otimes X^*
\stackrel{e_X \otimes id}{\longrightarrow}X^*$$
are equal to the identity morphisms.

(ii) A {\em left dual} to $X$ is an object $^*X$ with two morphisms
$$
e_X': X\otimes {}^*X\to \be,\; i_X': \be \to {}^*X\otimes X
$$
such that the compositions
$$ X\stackrel{i_X'\otimes id}{\longrightarrow}X\otimes {}^*X\otimes X
\stackrel{id \otimes e_X'}{\longrightarrow}X$$
$$ ^*X\stackrel{id\otimes i_X}{\longrightarrow}{}^*X\otimes X\otimes {}^*X
\stackrel{e_X' \otimes id}{\longrightarrow}{}^*X$$
are equal to the identity morphisms.

(iii) A monoidal category $\cC$ is called {\em rigid} if every object in 
$\cC$ has right and left duals.

{\bf Remark 2.} (see \cite{BK}) (i) Dual objects are defined canonically, that
is if there exists a dual (right or left), it is unique up to a unique 
isomorphism. 

(ii) For any object $X$ of rigid monoidal category $\cC$ there are canonical 
isomorphisms $X={}^*(X^*)=({}^*X)^*$.

(iii) Right (and left) duality can be canonically extended to a functor
$\cC \to \cC^{op}$ where $\cC^{op}$ is the opposite category to $\cC$. This
functor is equivalence of tensor categories.

Recall that we consider only semisimple categories. One can show that 
under this assumption for any object $X$ of a rigid monoidal category we
have a (noncanonical) isomorphism $^*X\simeq X^*$. Hence the Grothendieck
group $K_0(\cC )$ of a rigid monoidal category $\cC$ is a based ring.

If we assume that unit object of a rigid monoidal category $\cC$ is 
irreducible then $K_0(\cC)$ is unital based ring.

{\bf Conjecture 1.} (Ocneanu rigidity) {\em For a fixed finite dimensional 
unital based ring $R$ there are only finitely many rigid monoidal categories 
$\cC$ with $K_0(\cC )\simeq R$.}

As far as I know substantial progress in the proof of this Conjecture was
achieved by E.~Blanchard and A.~Wassermann, based on idea of A.~Ocneanu. One 
can also ask if a similar statement is true if we omit rigidity assumption.

{\bf Definition 5.} Let $\cC$ and $\cC'$ be two monoidal categories. A
{\em monoidal functor} $(F, b, u)$ is a triple consisting of a functor 
$F:\cC \to \cC'$, functorial isomorphism $b=\{ b_{X,Y}\},\; b_{X,Y}: 
F(X\otimes Y)=F(X)\otimes F(Y)$, and isomorphism $u: F(\be )=\be$ satisfying 
the natural compatibilities: the diagrams
$$
\xymatrix{F((X\otimes Y)\otimes Z)\ar[rr]^{b_{X\otimes Y,Z}}\ar[d]_{Fa_{X,Y,Z}}&&
F(X\otimes Y)\otimes F(Z)\ar[rr]^{b_{X,Y}\otimes id}&&(F(X)\otimes F(Y))\otimes
F(Z)\ar[d]^{a_{F(X),F(Y),F(Z)}}\\
F(X\otimes (Y\otimes Z))\ar[rr]^{b_{X,Y\otimes Z}}&&F(X)\otimes F(Y\otimes Z)
\ar[rr]^{id\otimes b_{Y,Z}}&&F(X)\otimes (F(Y)\otimes F(Z))}
$$
and
$$\xymatrix{F(\be \otimes X)\ar[r]^{b_{\be,X}}\ar[d]_{F(l_X)}&F(\be)\otimes 
F(X)\ar[d]^{u\otimes id}&&F(X\otimes \be)\ar[r]^{b_{X,\be}}\ar[d]_{F(r_X)}&
F(X)\otimes F(\be)\ar[d]^{id \otimes u}\\
F(X)&\be \otimes F(X)\ar[l]_{l_{F(X)}}&&F(X)&F(X)\otimes \be \ar[l]_{r_{F(X)}}}
$$
are commutative.

We give now several well known examples of monoidal categories.

{\bf Examples.} (i) The category $\Ve_k$ of finite dimensional vector spaces
over $k$ has a natural structure of a monoidal category where the functor 
$\otimes =\otimes_k$ is just the usual tensor product. This category is 
semisimple and rigid. The unit object is a one dimensional space $\be_k$ with 
fixed basis. The unit object is irreducible. For a monoidal category $\cC$ a 
{\em fiber functor} is a monoidal functor from $\cC$ to $\Ve_k$.

(ii) Let $G$ be an affine group scheme over $k$. Then category
$\Rep(G)$ of finite dimensional rational representations of $G$ has a natural
structure of a rigid monoidal category with irreducible unit object, which 
{\em is not} semisimple in general. The functor of forgetting the $G-$action 
has a natural structure of a fiber functor.

(iii) It is well known that the category of representations of a bialgebra is 
a monoidal category.
More generally, let $H$ be a weak bialgebra, see \cite{Sl, NV}. Then
the category $\Rep(H)$ of $H-$modules is a monoidal category.

(iv) Let $\cA$ be a semisimple abelian category. The category $\Fun(\cA, \cA)$
of functors from $\cA$ to $\cA$ has a structure of a monoidal category with
tensor product induced by composition of functors. This category is semisimple
and rigid (duality is given by taking adjoint functor). Its unit object is not
irreducible if $\cA$ has at least two nonisomorphic irreducible objects.

(v) If characteristic of the base field $k$ is not 2
there are exactly two categories with based ring isomorphic to 
$K_0(\Rep(\BZ/2\BZ))$, one is $\Rep(\BZ/2\BZ)$ itself and the second 
$\Rep(\BZ/2\BZ)^{tw}$ is new. In fact, in such a category there is only one
nontrivial associativity constraint (for triple product of nonunit object) and
in the category $\Rep(\BZ/2\BZ)^{tw}$ it differs by sign from the one in 
$\Rep(\BZ/2\BZ)$. Both categories are rigid. The category 
$\Rep(\BZ/2\BZ)^{tw}$ has no fiber functor. 

(vi) More generally, let $G$ be a finite group and consider the category
$\cC_G$ with (isomorphism classes of) simple objects $X_g$ parametrized by
$G$ and the tensor product functor given by $X_{g_1}\otimes X_{g_2}=
X_{g_1g_2}$. The monoidal structures on the category $\cC_G$ are parametrized
by the group $H^3(G,k^*)$, see e.g. \cite{Q,TY}.

\subsection{Module categories} The following definition is crucial for
this paper.

{\bf Definition 6.} A {\em module category} over a monoidal category $\cC$ is
a category $\cM$ together with an exact bifunctor $\otimes: \cC \times \cM \to 
\cM$ and functorial associativity and unit isomorphisms $m_{X,Y,M}:
(X\otimes Y)\otimes M
\to X\otimes (Y\otimes M),\; l_M: \be \otimes M\to M$ for any $X, Y\in \cC,\;
M\in \cM$ such that the diagrams
$$
\xymatrix{&((X\otimes Y)\otimes Z)\otimes M \ar[dl]_{a_{X,Y,Z}\otimes id} 
\ar[dr]^{m_{X\otimes Y,Z,M}}&\\
(X\otimes (Y\otimes Z))\otimes M \ar[d]^{m_{X,Y\otimes Z,M}}&&(X\otimes Y)
\otimes (Z\otimes M) \ar[d]_{m_{X,Y,Z\otimes M}}\\X\otimes ((Y\otimes Z)
\otimes M) \ar[rr]^{id\otimes m_{Y,Z,M}}&&X\otimes (Y\otimes (Z\otimes M))}
$$
and
$$
\xymatrix{(X\otimes \be )\otimes M\ar[rr]^{m_{X,\be ,Y}} \ar[dr]^{r_X\otimes 
id}&&X\otimes (\be \otimes M)
\ar[dl]_{id\otimes l_M}\\ &X\otimes M&}
$$
commute.

{\bf Remark 3.} As far as I know this definition first appeared in Bernstein's
lectures \cite{B} and in the work of L.~Crane and I.~Frenkel \cite{CF}. This
notion is implicitly present in Boundary Conformal Field Theory, see e.g.
\cite{BPPZ, FS, PZ}. In this context $\cC$ is the fusion category of the
corresponding Conformal Field Theory and irreducible objects of $\cM$ are
``boundary conditions''. It is clear that module categories can be described by
certain ``$6j-$symbols'' $^{(1)}F$ (this decription is analogous to the
$6j-$symbols description of monoidal categories). In Boundary Conformal Field
Theory these $6j-$symbols appear as coefficients of boundary field operator
product expansion. So we consider the notion of module category as a 
coordinate free version of Boundary Conformal Field Theory. Many examples of 
module categories (without using this name) were studied in Operator Algebras 
Theory, see e.g. \cite{BEK, O}.

The Grothendieck group $K_0(\cM)$ of a module category $\cM$ over a monoidal
category $\cC$ with basis given by classes of irreducible objects is clearly
a $\BZ_+-$module over the $\BZ_+-$ring $K_0(\cC)$.

{\bf Lemma 2.} Let $\cC$ be a rigid monoidal category and $\cM$ be a module
category over $\cC$. Then for any $X\in \cC,\; M_1, M_2\in \cM$ we have 
canonical isomorphisms
$$\Hom(X\otimes M_1, M_2)\cong \Hom(M_1, {}^*X\otimes M_2),\;
\Hom(M_1, X\otimes M_2)=\Hom(X^*\otimes M_1,M_2).$$

{\bf Proof.} Clear. \sq

The Lemma implies that for a module category $\cM$ over a rigid monoidal 
category $\cC$ the Grothendieck group $K_0(\cM)$ is a based module over the 
based ring $K_0(\cC)$.

{\bf Examples.} (i) Any monoidal category $\cC$ has a structure of a module
category over itself with associativity and unit isomorphisms induced by the
ones from the monoidal structure on $\cC$.

(ii) Let $F: \cC \to \Ve_k$ be a fiber functor. It defines a structure of a
module category over $\cC$ on the category $\Ve_k$ as follows: for $X\in \cC$
and $V\in \Ve_k$ we set $X\otimes V:=F(X)\otimes_k V$ with associativity 
and unit isomorphisms defined as compositions:
$$ (X\otimes Y)\otimes V=F(X\otimes Y)\otimes_k V\stackrel{b_{X,Y}\otimes_k id}
{\longrightarrow}(F(X)\otimes_k F(Y))\otimes_k V=$$
$$
=F(X)\otimes_k (F(Y)\otimes_k 
V)=X\otimes (Y\otimes V)$$
and
$$\be \otimes V=F(\be)\otimes_k V\stackrel{u\otimes_k id}{\longrightarrow}
\be_k\otimes V=V.$$
It is easy to see that conversely, a structure of module category over $\cC$ 
on $\Ve_k$ determines a fiber functor on $\cC$ (see Section 4 for a more 
general statement).

(iii) Assume for a moment that the field $k$ is not algebraically
closed and consider the category $\cC=\Ve_k$. It is easy to see that 
indecomposable module categories over $\cC$ are classified by (finite
dimensional) skew fields over $k$, or equivalently by the Brauer groups
of all finite extensions of $k$. 

{\bf Definition 7.} (i) Let $\cM_1$ and $\cM_2$ be two module categories over
a monoidal category $\cC$. A {\em module functor} from $\cM_1$ to $\cM_2$
is a functor $F: \cM_1\to \cM_2$ together with functorial morphism
$c_{X,M}: F(X\otimes M)\to X\otimes F(M)$ for any $X\in \cC,\; M\in \cM_1$
such that the diagrams
$$
\xymatrix{&F((X\otimes Y)\otimes M)\ar[dl]_{Fm_{X,Y,M}} \ar[dr]^{c_{X\otimes 
Y}}&\\ F(X\otimes (Y\otimes M))\ar[d]^{c_{X,Y\otimes M}}&&(X\otimes Y)
\otimes F(M)\ar[d]_{m_{X,Y,F(M)}}\\X\otimes F(Y\otimes M)\ar[rr]^{id\otimes 
c_{Y,M}}&&X\otimes (Y\otimes F(M)}
$$
and
$$
\xymatrix{F(\be \otimes M)\ar[rr]^{Fl_M} \ar[dr]^{c_{\be,M}}&&F(M)\\ 
&\be \otimes F(M)\ar[ur]^{l_{F(M)}}&}
$$
are commutative.

(ii) We say that a module functor $(F, c_{X,M})$ is {\em strict} if all
morphisms $c_{X,M}$ are isomorphisms.

(iii) Two module categories $\cM_1$ and $\cM_2$ over $\cC$ are {\em equivalent}
if there exists a strict module functor from $\cM_1$ to $\cM_2$ which is an
equivalence of categories.

(iv) For two module categories $\cM_1$ and $\cM_2$ over a monoidal category
$\cC$ their {\em direct sum} is the category $\cM_1\times \cM_2$ with 
coordinatewise additive and module structure.

(v) A module category is {\em indecomposable} if it is not equivalent to a
direct sum of two non trivial module categories.

{\bf Remark 4.} For a rigid monoidal category $\cC$ any module functor is
automatically strict.

{\bf Conjecture 2.} {\em For a given rigid monoidal category with finitely many
irreducible objects there exists only finitely many inequivalent 
indecomposable module categories.}

\section{Morita theory for module categories}
In this section $\cC$ denotes a semisimple monoidal category.
\subsection{Algebras in monoidal categories}\par

{\bf Definition 8.} (cf. e.g. \cite{Br}) (i) An {\em algebra} in a monoidal 
category $\cC$ is an object $A$ of $\cC$ endowed with 
a {\em multiplication morphism} 
$m:A\otimes A\to A$ and a {\em unit morphism} $e: \be \to A$ such that the 
diagrams
$$\xymatrix{&(A\otimes A)\otimes A\ar[rd]^{m\otimes id}\ar[ld]_{a_{A,A,A}}&\\
A\otimes (A\otimes A)\ar[d]_{id\otimes m}&&A\otimes A\ar[d]^m\\
A\otimes A\ar[rr]^m&&A}$$
and
$$\xymatrix{&\be \otimes A\ar[dl]_{l_A}\ar[dr]^{e\otimes id}&\\
A&&A\otimes A\ar[ll]_m}
\xymatrix{&A \otimes \be\ar[dl]_{r_A}\ar[dr]^{id\otimes e}&\\
A&&A\otimes A\ar[ll]_m}$$
commute.

(ii) A {\em right module} over an algebra $A$ in a monoidal category $\cC$ is
an object $M$ of $\cC$ together with an {\em action morphism} 
$a: M\otimes A\to M$ such that the diagrams
$$\xymatrix{M\otimes A\otimes A\ar[r]^{id\otimes m}\ar[d]_{a\otimes id}&
M\otimes A\ar[d]^a\\ M\otimes A\ar[r]^a&M}
\xymatrix{&M \otimes \be\ar[dl]_{r_M}\ar[dr]^{id\otimes e}&\\
M&&M\otimes A\ar[ll]_a}$$
commute. A {\em left module} over $A$ is defined in a similar way.

(iii) A {\em morphism} between two right modules $M_1$ and $M_2$ over $A$ is
a morphism $\alpha \in \Hom_{\cC}(M_1,M_2)$ such that the diagram
$$\xymatrix{M_1\otimes A\ar[r]^{\alpha \otimes id}\ar[d]_{a_1}&
M_2\otimes A\ar[d]^{a_2}\\M_1\ar[r]^\alpha &M_2}$$
commutes. It is clear that morphisms between $M_1$ and $M_2$ form a
$k-$subspace of $\Hom_{\cC}(M_1,M_2)$. We will denote this subspace
$\Hom_A(M_1,M_2)$.

{\bf Exercise.} Prove that if $M$ is a right $A-$module then $^*M$ (but not
$M^*$!) has a natural structure of a left $A-$module. If $M$ is a left 
$A-$module then $M^*$ has natural structure of a right $A-$module. These two 
functors are inverse to each other.

We leave to the reader to define morphisms between algebras, ideals, bimodules,
tensor products over $A$ etc.

{\bf Lemma 3.} Let $A$ be an algebra in a monoidal category $\cC$. Right (and
similarly left) modules over $A$ with morphisms defined above form an abelian
category $\Mod_{\cC}(A)$.

{\bf Proof.} Clear. \sq

For any right module $M$ over an algebra $A$ in $\cC$ and object $X$ of $\cC$
the object $X\otimes M$ has natural structure of $A-$module with action
morphism given by $id\otimes a$. Moreover, for any $X, Y\in \cC$ associativity
isomorphism between $(X\otimes Y)\otimes M$ and $X\otimes (Y\otimes M)$ is
a morphism of $A-$modules. It is straightforward to check that in this way
we endow $\Mod_{\cC}(A)$ with a structure of module category over $\cC$.
We will prove in this section that if $\cC$ is rigid and has finitely many
irreducible objects then any semisimple module category is equivalent to 
$\Mod_{\cC}(A)$ for some algebra $A$.

{\bf Definition 9.} (i) An algebra $A$ in $\cC$ is {\em semisimple} if the 
category $\Mod_{\cC}(A)$ is semisimple.

(ii) An algebra $A$ is called {\em indecomposable} if the module category 
$\Mod_{\cC}(A)$ is indecomposable.

{\bf Remark 5.} It is easy to see that $A$ is indecomposable in our sense if
and only if it is indecomposable in usual sense --- it is not isomorphic to
a direct sum of nontrivial algebras.

{\bf Definition 10.} Two algebras $A_1$ and $A_2$ in $\cC$ are {\em Morita
equivalent} if module categories $\Mod_{\cC}(A_1)$ and $\Mod_{\cC}(A_2)$ are
equivalent. 

{\bf Remark 6.} Another definition (not using the notion of module categories)
of Morita equivalence was proposed in \cite{FS}. It is easy to see that our 
definition and the definition of \cite{FS} are the same. We hope that our
language gives a little bit more flexibility (see e.g. below the proof
of Theorem 6).

Let $A$ be an algebra in $\cC$. For any $X\in \cC$ the object $X\otimes A$
has a natural structure of right $A-$module induced by multiplication in $A$.

{\bf Lemma 4.} For any $A-$module $L$ and an object $X\in \cC$ we have
canonical isomorphism
$$\Hom_A(X\otimes A,L)=\Hom(X,L).$$

{\bf Proof.} It is clear that unit morphism $e:\be \to A$ defines canonical
isomorphism $\Hom_A(A,L)=\Hom(\be,L)$. Now the isomorphism of Lemma can be
obtained as composition
$$\Hom_A(X\otimes A,L)=\Hom_A(A,{}^*X\otimes L)=\Hom(\be,{}^*X\otimes L)=
\Hom(X,L).$$ \sq

{\bf Remark 7.} In the proof above we used the rigidity of $\cC$. In fact this 
can be avoided by a more lengthy calculation. We leave this to the reader as 
an exercise.

We see that modules $X\otimes A$ are projective objects of $\Mod_{\cC}(A)$ and
any irreducible $A-$module is a quotient of a module of the form $X\otimes A$.
In particular the category $\Mod_{\cC}(A)$ has enough projective objects.

\subsection{Internal Hom for module categories} In this subsection $\cC$ is a
semisimple rigid monoidal category and $\cM$ is a semisimple module category
over $\cC$.

{\bf Definition 11.} Let $M_1$ and $M_2$ be two objects of $\cM$. Their
{\em internal} Hom is an ind-object $\iHom(M_1,M_2)$ of $\cC$ representing
functor $X\mapsto \Hom(X\otimes M_1,M_2)$. 

{\bf Remark 8.} (i) Functor $X\mapsto \Hom(X\otimes M_1,M_2)$ is exact whence
existence of $\iHom(M_1,M_2)$.

(ii) If both categories $\cC$ and $\cM$ have finitely many irreducible
objects, then $\iHom(M_1,M_2)$ is an object of $\cC$.

(iii) By Yoneda's Lemma the object $\iHom(M_1,M_2)$ is uniquely defined up
to unique isomorphism, so $\iHom(?,?)$ is a bifunctor.

{\bf Lemma 5.} We have canonical isomorphisms
$$\begin{array}{cl}
(1)&\Hom(X\otimes M_1,M_2)=\Hom(X,\iHom(M_1,M_2)),\\
(2)&\Hom(M_1,X\otimes M_2)=\Hom(\be,X\otimes \iHom(M_1,M_2)),\\
(3)&\iHom(X\otimes M_1,M_2)=\iHom(M_1,M_2)\otimes X^*,\\
(4)&\iHom(M_1,X\otimes M_2)=X\otimes \iHom(M_1,M_2).\end{array}$$

{\bf Proof.} Formula (1) is just the definition of $\iHom(M_1,M_2)$ and 
isomorphism (2) is the composition
$$\Hom(M_1,X\otimes M_2)\cong \Hom(X^*\otimes M_1,M_2)=$$
$$=\Hom(X^*,\iHom(M_1,M_2))\cong \Hom(\be, X\otimes \iHom(M_1,M_2)).$$
We get isomorphism (3) from the calculation
$$\Hom(Y,\iHom(X\otimes M_1,M_2))=\Hom(Y\otimes (X\otimes M_1),M_2)=
\Hom((Y\otimes X)\otimes M_1,M_2)=$$
$$=\Hom(Y\otimes X,\iHom(M_1,M_2))=\Hom(Y,\iHom(M_1,M_2)\otimes X^*)$$
and isomorphism (4) from the calculation
$$\Hom(Y,\iHom(M_1,X\otimes M_2))=\Hom(Y\otimes M_1,X\otimes M_2)=$$
$$=\Hom(X^*
\otimes (Y\otimes M_1),M_2)=\Hom((X^*\otimes Y)\otimes M_1,M_2)=$$
$$=\Hom(X^*\otimes Y,\iHom(M_1,M_2))=\Hom(Y,X\otimes \iHom(M_1,M_2)).$$ \sq

For two objects $M_1, M_2$ of $\cM$ we have the canonical morphism 
$$ev_{M_1,M_2}: \iHom(M_1,M_2)\otimes M_1\to M_2$$ 
obtained as the image of $id$ under the isomorphism
$$\Hom(\iHom(M_1,M_2),\iHom(M_1,M_2))=\Hom(\iHom(M_1,M_2)\otimes M_1,M_2).$$
Let $M_1, M_2, M_3$ be three objects of $\cM$. Then there is a canonical
composition morphism
$$(\iHom(M_2,M_3)\otimes \iHom(M_1,M_2))\otimes M_1=\iHom(M_2,M_3)\otimes
(\iHom(M_1,M_2)\otimes M_1)\stackrel{id\otimes ev_{M_1,M_2}}{\longrightarrow}$$
$$\stackrel{id\otimes ev_{M_1,M_2}}{\longrightarrow}
\iHom(M_2,M_3)\otimes M_2\stackrel{ev_{M_2,M_3}}{\longrightarrow}M_3$$
which produces the {\em multipication morphism} 
$$\iHom(M_2,M_3)\otimes \iHom(M_1,M_2)\to \iHom(M_1,M_3)$$ 
(note that order of factors is opposite to 
intuitive one!). It is straightforward to check that this multiplication is
associative and compatible with the isomorphisms of Lemma 5.

\subsection{Main Theorem} In this subsection we assume that $\cC$ is a 
semisimple rigid monoidal category with finitely many irreducble objects and 
irreducible unit object.

{\bf Theorem 1.} Let $\cM$ be a semisimple indecomposable module category over
$\cC$. Then there exists a semisimple indecomposable algebra $A\in \cC$ such 
that the module categories $\cM$ and $\Mod_{\cC}(A)$ are equivalent.

{\bf Proof.} Fix any nonzero object $M$ of $\cM$. The multiplication morphism
defines a structure of an algebra on $A=\iHom(M,M)$. Consider a functor $F$ 
from $\cM$ to $\cC$ defined by $N\mapsto \iHom(M,N)$. Again the multiplication
morphism defines a structure of a right $A-$module on $\iHom(M,N)$ and hence
$F$ is a functor from $\cM$ to $\Mod_{\cC}(A)$. Isomorphism (4) of
Lemma 5 defines a structure of a module functor on $F$ (axioms of a module
functor follow from compatibility of the multiplication with isomorphism (4) 
of Lemma 5).

Now we claim that the functor $F:\cM \to \Mod_{\cC}(A)$ is an equivalence of 
categories. We will proceed in steps:

(1) If $N\ne 0$ then $F(N)\ne 0$. 

Indeed, otherwise category $\cM$ is decomposable (objects $N\in \cM$ such
that $F(N)=0$ clearly form a module subcategory which is a module direct
summand of $\cM$ thanks to the rigidity of $\cC$, see Lemma 1).

(2) The functor $F$ is injective on Hom's.

This follows immediately from the semisimplicity $\cM$ and (1).

(3) The map $F: \Hom(N_1,N_2)\to \Hom_A(F(N_1),F(N_2))$ is surjective
(and hence isomorphism by (2)) for any $N_2\in \cM$ and $N_1$ of the form
$X\otimes M,\; X\in \cC$.

Indeed, $F(N_1)=\iHom(M,X\otimes M)=X\otimes A$ and the statement follows from
the calculation:
$$\Hom_A(F(N_1),F(N_2))=\Hom_A(X\otimes A,F(N_2))=\Hom(X,F(N_2))=$$
$$=\Hom(X,\iHom(M,N_2))=\Hom(X\otimes M,N_2)=\Hom(N_1,N_2).$$

(4) The map $F: \Hom(N_1,N_2)\to \Hom_A(F(N_1),F(N_2))$ is an isomorphism for
any $N_1, N_2\in \cM$.

It is clear that there exist objects $X,Y\in \cC$ and an exact sequence
$$Y\otimes M\to X\otimes M\to N_1\to 0.$$
Hence (4) is consequence of (3).

(5) The functor $F$ is surjective on isomorphism classes of objects of 
$\Mod_{\cC}(A)$.

We know from Lemma 4 that for any object $L\in \Mod_{\cC}(A)$ there
exists an exact sequence
$$Y\otimes A\stackrel{\tilde f}{\longrightarrow} X\otimes A\to L\to 0$$
for some $X, Y\in \cC$. Let $f\in \Hom(Y\otimes M,X\otimes M)$ be the preimage
of $f$ under the isomorphism
$$\Hom(Y\otimes M,X\otimes M)=\Hom_A(F(Y\otimes M),F(X\otimes M))=\Hom_A
(Y\otimes A,X\otimes A)$$
and let $N\in \cM$ be the cokernel of $f$. It is clear that $F(N)=L$.

We proved that $F$ is equivalence of categories and proved the Theorem. \sq

{\bf Remark 9.} (i) If a category $\cC$ has infinitely many simple objects,
one can prove a similar Theorem by considering ind-algebras and ind-modules.

(ii) The proof of the fact that $F$ is an equivalence of categories follows
the standard pattern from homological algebra.

{\bf Examples.} (i) If $\cM=\cC$ is the ``regular'' module category
and $X\in \cC$ it is easy to see that $\iHom(X,X)=X\otimes X^*$.

(ii) Let $\cM =\Ve_k$ and thus is associated with fiber functor $F:\cC \to 
\Ve_k$. By the usual Tannakian formalism this induces an equivalence 
$\cC =\Rep(H)$
for some Hopf algebra $H$. In this case $\iHom(k,k)=H^*$ --- the dual Hopf
algebra together with the natural $H-$action.

{\bf Corollary} (of proof). Any semisimple indecomposable algebra in $\cC$ is 
Morita equivalent to the algebra $A$ with $\Hom(\be,A)=k$.

{\bf Proof.} Indeed it is enough to take for $M$ a simple object of $\cM$. \sq

\subsection{Example: module categories over $\Rep(G)$} Let $G$ be a finite
group. Assume that the base field $k$ is algebraically closed and the
characteristic of $k$ does not divide order of $G$, so
the category $\cC =\Rep(G)$ is semisimple. In this subsection we classify all
module categories over $\Rep(G)$.

{\bf Example.} Let $H\subset G$ be a subgroup, let $1\to k^*\to \tilde H\to 
H\to 1$ be a central extension of $H$ whose kernel is identified with 
the multiplicative group $k^*$. Then the category $\Rep^1(\tilde H)$ of 
representations $V$ of $\tilde H$ such that $k^*$ acts on $V$ via identity 
character is a module category over $Rep(G)$ via usual tensor product. The
extensions $\tilde H$ as above are in one to one correspondence with elements
of the group $H^2(H, k^*)$. For an element $\omega \in H^2(H, k^*)$ we will
write $\Rep^1(H,\omega)$ instead of $\Rep^1(\tilde H)$ where $\tilde H$
corresponds to $\omega$.

{\bf Theorem 2.} (cf. \cite{BO}) Any indecomposable module category over 
$\Rep(G)$ is equivalent to $\Rep^1(H,\omega)$ for some $H\subset G$ and
$\omega \in H^2(H, k^*)$. Two module categories $\Rep(H_1,\omega_1)$ and
$\Rep(H_2,\omega_2)$ are equivalent if and only if pairs $(H_1,\omega_1)$
and $(H_2,\omega_2)$ are conjugate under the adjoint action of $G$.

{\bf Proof.} First it is easy to see that for any $V\in \Rep^1(\tilde H)$
we have $\iHom(V,V)=\Ind_H^G\End(V)$.

Now let $\cM$ be a semisimple indecomposable module category over $\Rep(G)$.
By Theorem 1 $\cM$ is equivalent to $\Mod_{\Rep(G)}(A)$ for an indecomposable
semisimple $G-$algebra $A$. It is easy to see that semisimplicity of a
$G-$algebra $A$ implies semisimplicity of $A$ as an algebra in $\Ve_k$, see
\cite{BO}. So $A$ is just a direct sum of matrix algebras. The group $G$ acts
on the set of minimal central idempotents of $A$ and this action is transitive
since $A$ is indecomposable. Let $e$ be a minimal central idempotent and let
$H$ be its stabilizer in $G$. It is clear that the subalgebra $eAe$ is 
$H-$invariant and $A=\Ind_H^G(eAe)$. The algebra $eAe$ is a matrix algebra and
hence $eAe=\End(V)$ for some projective representation of $H$ since all
automorphisms of a matrix algebra are inner. The Theorem is proved. \sq

{\bf Examples.} (i) Let $G=\BZ /2\BZ \times \BZ /2\BZ$. The only subgroup of 
$G$ having central extensions is $G$ itself. The group $G$ has exactly one
irreducible projective representation (of dimension 2) even though it has
two different central extensions. We see that the category $\Rep(G)$ has two
different module categories with one irreducible object, hence the category
$\Rep(G)$ has an additional fiber functor. In general we see that fiber 
functors $\Rep(G)\to \Ve$ are classified by conjugacy classes of pairs
$(H, \omega)$ such that the category $\Rep^1(H,\omega)$ has only one 
irreducible object, and in particular the order of $H$ is a square. This 
result is due to M.~V.~Movshev, see \cite{Mo}.

(ii) {\bf Definition 12.} (\cite{EG}) Two finite groups
$G_1$, $G_2$ are called {\em isocategorical} if $\Rep(G_1)$ is equivalent to
$\Rep(G_2)$ as a monoidal category.

Of course if the groups $G_1$, $G_2$ are isocategorical then the Grothendieck 
rings $K_0(\Rep(G_1))$ and $K_0(\Rep(G_2))$ are isomorphic as based rings (or
equivalently the character tables of $G_1$ and $G_2$ are the same). But the 
property of being isocategorical is much stronger. For example, it is known
(see \cite{TY}, \cite{EG}) that	the two nonabelian groups of order 8 are not
isocategorical (here is a simple proof of this fact: let us calculate the
number of fiber functors for these categories; we always have the tautological
fiber functor, and it is easy to see from the above that the additional fiber
functors are classified by conjugacy classes of subgroups isomorphic to
$\BZ/2\BZ \times \BZ/2\BZ$. For the quaternion group we have no such subgroups
and for the dihedral group we have two conjugacy classes of such subgroups.
So the corresponding monoidal categories have a different number of fiber
functors and are not equivalent). 

Let $G$ be a finite group. Let $\cM(G)$ be a finite set consisting of pairs
$(H,\omega)$ where $H\subset G$ is a subgroup and $\omega \in H^2(H,k^*)$ is
a cohomology class. Let $m:\cM(G)\to \BN$ be a function given by $m((H,\omega))
:=\mbox{number of irreducible objects in}\; \Rep^1(H,\omega)$. It follows from
the classification of module categories over $\Rep(G)$ above that for two
isocategorical groups $G_1$ and $G_2$ there is a bijection $\cM(G_1)\to
\cM(G_2)$ preserving the function $m$. 
In \cite{EG} P.~Etingof and S.~Gelaki described all pairs of nonisomorphis
isocategorical groups and in particular constructed explicit examples. It
would be interesting to describe the bijection above for these examples.

\subsection{Semisimplicity} In this section we discuss the somewhat subtle 
question of semisimplicity of algebra $A$.

{\bf Proposition 2.} (i) A splitting of the multiplication morphism
$m: A\otimes A\to A$ as a morphism of bimodules over $A$ is sufficient for 
semisimplicity of algebra $A$.

(ii) Assume that $\Hom(\be,A)=k$ and let $\varepsilon \in \Hom(A,\be)$ be a
nonzero morphism. Semisimplicity of $A$ implies that pairing
$$A\otimes A\stackrel{m}{\longrightarrow}A\stackrel{\varepsilon}
{\longrightarrow}\be$$
is nondegenerate; that is, this map defines an isomorphism $A\to A^*$.

{\bf Proof.} (i) Indeed, this condition implies that any $A-$module $L=
L\otimes_AA$ is a direct summand of $L\otimes A=L\otimes_A(A\otimes A)$ and
consequently is projective.

(ii) It is clear that the map $A\to A^*$ is a morphism of right $A-$modules.
In the semisimple case the (right or left) $A-$module $A$ is irreducible since
$\Hom_A(A,A)=\Hom(\be,A)=k$ and the map $A\to A^*$ is isomorphism by the Schur
Lemma. \sq

{\bf Example.} Condition (i) is not necessary as the following example shows. 
Let $k$ be of characteristic 2 and let $H$ be a Hopf algebra dual to the group
algebra of $G=\BZ /2\BZ$. It is easy to see that $k[G]$ is semisimple as 
an $H-$algebra but does not satisfy condition (i). 

It is likely that the condition of nondegeneracy from (ii) is not sufficient 
for semisimlicity in general, but I don't know counterexample.

Suppose now that in a category $\cC$ we have functorial isomorphisms
$\delta_X:X\to X^{**}$ satisfying
$$\delta_{X\otimes Y}=\delta_X\otimes \delta_Y,\; \delta_{\be}=id,\; 
\delta_{X^*}=(\delta_X^*)^{-1} $$
where for $f\in \Hom(X,Y), \;f^*\in \Hom(Y^*,X^*)$ is transposed morphism.
In this case we can identify $X^*$ and $^*X$, and quantum dimension is
defined. The following Theorem is proved in \cite{KO} Theorem 4.3 (we work 
there in a braided category $\cC$ with a commutative algebra $A$, but the 
proof of 4.3 does not use neither brading nor commutativity of $A$):

{\bf Theorem 3.} Assume that $\cC$ is as above. Let $A$ be an algebra in $\cC$
with $\Hom(\be,A)=k$ and with a nongenerate pairing defined above. Assume in
addition that $\dim(A)\ne 0$. Then condition (i) of Proposition 2 holds and
hence the algebra $A$ is semisimple. \sq

{\bf Remark 10.} (i) For monoidal categories coming from Conformal Field Theory
condition $\dim(A)\ne 0$ is automatically satisfied since all dimensions are
positive.

(ii) Being motivated by the Boundary Conformal Field Theory J.~Fuchs and 
C.~Schweigert introduced the notion of a {\em Frobenius algebra}, see
\cite{FS}. By definition a Frobenius algebra $A$ is an algebra together with 
a splitting morphism $\Delta: A\to A\otimes A$ which is assumed in addition to
be coassociative. Note that the splitting constructed in Theorem 3 is 
automatically coassociative since it is just the dualization of the associative
map $m: A\otimes A\to A$.

(iii) We sketch here a possible way to approach Conjecture 2 for monoidal
categrories coming from CFT. In view of Lemma 1 and Proposition 1 it would
be enough to show that module categories have no deformations. In view of
Theorem 1 this reduces to showing that semisimple algebras have no
deformations. Such infinitisemal deformations should be described by
cohomology of Hochschild complex of $A-A$ bimodule $A$. But by Proposition 2
and Theorem 3 above this bimodule is projective, hence its cohomology vanishes.

\section{Weak Hopf algebras}
\subsection{} Let $H$ be a bialgebra. The category $\Rep(H)$ has a natural 
monoidal structure, and many examples of monoidal categories arise in this way.
But it is well known that not any monoidal category is equivalent to $\Rep(H)$
for some bialgebra $H$. The reason is the following: the forgetful functor 
$\Rep(H)\to \Ve_k$ is clearly a fiber functor but for a general monoidal 
category there is no reason to have a fiber functor. Next thing to try is the 
following: let $R$ be a separable algebra (that is $R$ is a finite direct sum 
of matrix algebras) and consider the category of finite dimensional 
$R-$bimodules $\Bimod(R)$. The category $\Bimod(R)$ is a monoidal category 
with monoidal structure induced by the tensor product over $R$.

{\bf Defintion 13.} An {\em $R-$fiber functor} for a monoidal category $\cC$ 
is a monoidal functor $\cC \to \Bimod(R)$.

Let $\cC$ be a monoidal category and let $F:\cC \to \Bimod(R)$ be a $R-$fiber
functor. It is natural to expect that Tannakian formalism works in such
situation and the functor $F$ induces an equivalence $\cC \to \Rep(H)$, where 
$H$ is some generalization of a bialgebra. This is indeed true and the
corresponding structure on $H$ is a structure of a {\em weak bialgebra}, see
\cite{Sl}. Recall here the definition of a weak bialgebra.

{\bf Definition 14.} (see \cite{NV, N, Sl}) A weak bialgebra is a finite
dimensional vector space $H$ with the structures of an associative algebra 
$(H, m, 1)$ with multiplication $m: H\otimes H\to H$ and unit $1\in H$ and a 
coassociative coalgebra $(H, \Delta, \eps)$ with comultiplication $\Delta: 
H\to H\otimes H$ and counit $\eps: H\to k$ such that: 

(i) The comultiplication $\Delta$ is a homomorphism of algebras such that
$$(\Delta \otimes id)\Delta (1)=(\Delta(1)\otimes 1)(1\otimes \Delta (1))=
(1\otimes \Delta (1))(\Delta (1)\otimes 1),$$

(ii) The counit satisfies the identity:
$$\eps(fgh)=\eps(fg^{(1)})\eps(g^{(2)}h)=\eps(fg^{(2)})\eps(g^{(1)}h)$$
for all $f, g, h\in H$.

The distinction between the definitions of a bialgebra and of a weak
bialgebra is the following: in the definition of weak bialgebra it is not
assumed that the coproduct preserves the unit and dually it is not assumed that
the counit is an algebra homomorphism. The algebra $R$ (denoted by $H_t$ in
\cite{NV}) is called the {\em base} of a weak bialgebra $H$. We refer the 
reader to \cite{NV} and \cite{Sl} for a review of the theory of weak Hopf 
algebras (which are weak bialgebras with an antipode). A relation between the 
theory of weak Hopf algebras and module categories is given by the following 

{\bf Proposition 3.} Let $\cC$ be a monoidal category and let $R$ be a 
separable algebra. There is natural bijection between sets \{ $R-$fiber 
functors\} and \{ structures of a module category over $\cC$ on $\Rep(R)$\} .

{\bf Proof.} The category $\Rep(R)$ has an obvious structure of a module 
category over $\Bimod(R)$ where the bifunctor $\Bimod(R)\otimes \Rep(R)\to 
\Rep(R)$ is the tensor product over $R$. Hence any $R-$fiber functor induces a
structure of a module category on $\Rep(R)$.

On the other hand, let us reformulate the definition of a module category in 
the following way: For an abelian category $\cM$, let $\Fun(\cM,\cM)$ denote 
the category of exact functors $\cM\to \cM$ with natural transformations as 
morphisms. The category $\Fun(\cM,\cM)$ has a monoidal structure with tensor 
product being the composition of functors. Now suppose that $\cM$ is a module
category over a monoidal category $\cC$. Any object $X\in \cC$ defines a 
functor $F_X: \cM \to \cM,\; F_X(M)=X\otimes M$. So we have a functor 
$F:\cC \to \Fun(\cM,\cM)$. The associativity isomorphism defines a natural 
transformation of functors $F_{X\otimes Y}\to F_X\circ F_Y$, and one checks 
that the axioms of a module category are equivalent to saying that $F$ is a 
monoidal functor. Now structures of module category over $\cC$ on $\Rep(R)$ 
are the same as monoidal functors $\cC \to \Fun(\Rep(R),\Rep(R))=\Bimod(R)$.

We leave it to the reader to check that the two constructions above are 
mutually inverse. \sq 

As an easy application we get the following statement:

{\bf Theorem 4.} Let $\cC$ be a semisimple monoidal category with finitely
many simple objects. Then there exists an equivalence $\cC \cong \Rep(H)$
where $H$ is a weak bialgebra. Moreover, we can assume that the base of $H$ is
commutative.

{\bf Proof.} As we mentioned above, $\cC$ is a module category over itself.
Choose an algebra $R$ such that $\cC$ is equivalent to $\Rep(R)$ as abelian
category. By the Proposition above we get $R-$fiber functor $\cC \to 
\Bimod(R)$. By the Theorem 1.8 from \cite{Sl} (see also section 1.3 there) we
get an equivalence $\cC \to \Rep(H)$ where $H$ is a weak bialgebra. \sq

{\bf Remark 11.} (i) As it was pointed to me by L.~Vainerman, this Theorem
was proved previously by T.~Hayashi, see \cite{H}.

(ii) Some version of the Theorem 4 is known in Operator Algebras
theory and in physics, see e.g. \cite{O, BEK, PZ}. Note that in \cite{O, BEK}
different terminology is used: there weak Hopf algebras are replaced by
closely related objects --- double triangular algebras, see \cite{PZ}.

(iii) Note that a weak bialgebra constructed in Theorem 4 is somewhat
noncanonical; the canonical object is a bialgebroid over $R$ (see \cite{Sl}). 
To get a weak bialgebra from a bialgebroid, one has to choose a separability 
idempotent in $R$, see \cite{Sl}. If $R$ is abelian then a separability
idempotent is unique (and our weak bialgebra is canonical); moreover over a
field of characteristic 0 the canonical choice of a separability idempotent
is possible since a symmetric separability idempotent is unique.

(iv) If the category $\cC$ is rigid, one easily defines the
antipode map $S: H\to H$. It is proved in \cite{N} that the map $S$
satisifies the axioms of the antipode in a weak Hopf algebra and so
$H$ becomes the weak Hopf algebra. 

\subsection{Duality for weak Hopf algebras} \label{dual}
The definition of a weak bialgebra has a virtue of being selfdual, that is if
$H$ is a weak bialgebra then so is $H^*$, see \cite{NV, Sl}. In this section
we explain the categorical meaning of this duality. Let $\cM =\Rep(R)$ be a
module category over a monoidal category $\cC$. Let $H$ be the corresponding
weak bialgebra constructed by Proposition 3. Consider the category $\cC^*:=
\Fun_\cC (\cM,\cM)$ of module functors from $\cM$ to itself. The category
$\cC^*$ is an abelian (non-semisimple, in general) monoidal category with 
composition of functors as a tensor product. Evidently, $\cM$ is a 
module category over $\cC^*$.

{\bf Theorem 5.} There is a natural monoidal equivalence of categories $\cC^*$ 
and $\Rep(H^*)$.

{\bf Proof.} Any functor $\Rep(R)\to \Rep(R)$ is isomorphic to the functor
$V\otimes_R?$ for some $R-$bimodule $V$. By definition an object of $\cC^*$
is an $R-$bimodule $V$ together with a functorial isomorphism $c_{X,M}:
V\otimes_R(X\otimes_RM)\to X\otimes_R(V\otimes_RM)$ for any $X\in \cC, M\in
\Rep(R)$. Functoriality in the variable $M$ implies that these isomorphisms
are induced by $R-$bimodule isomorphisms $c_X: V\otimes_RX\to X\otimes_RV$.
The isomorphisms $c_X$ should satisfy two conditions: $c_{X_1\otimes X_2}=
(id\otimes c_{X_2})(c_{X_1}\otimes id)$ (the pentagon diagram in Definition 7)
and $c_{\be}=id$ (the triangle axiom in Definition 7). The functotiality in
the variable $X$ implies that the isomorphisms $c_X$ are completely defined by 
the isomorphism $c_H$ (where $H$ is considered as an object of $\cC=\Rep(H)$)
and moreover by the restriction $c: V\to H\otimes_RV$ of this isomorphism to 
$V\otimes 1$ since any vector $x\in X$ is an image of $1\in H$ under a unique
$H-$module morphism $H\to X$. One verifies easily that the conditions on
the isomorphisms $c_X$ above are equivalent to the condition that $c$ defines
a structure of an $H-$comodule on $V$ (in fact this is only a comodule in the
category $\Bimod(R)$, to get a genuine comodule one uses a separability
idempotent in $R$ to imbed $H\otimes_RV\subset H\otimes V$). Conversely, let
$V$ be an $H-$comodule ($=$ module over $H^*$). For any object $X$ of $\cC$
define a map $c_X: V\otimes_RX\to X\otimes_RV$ by the formula 
$$c_X(v\otimes x)=\sum_ih_ix\otimes v_i$$
where $v\mapsto \sum_ih_i\otimes v_i$ is a coaction of $H$ on the element
$v\in V$. One verifies immediately that this map satisfies the conditions
above if it is well defined. So the proof of Theorem 5 is completed modulo
the following 

{\bf Lemma 6.} (D.~Nikshych) The map $c_X$ above is a well defined map of 
$R-$bimodules.

{\bf Proof.} We will use in the proof the notations from \cite{NV}. 

(1) The map $c_X$ is well defined: let $z\in H_t$ be an arbitrary element.
It acts on $v\in V$ by the formula $v\mapsto \sum_i\eps(h_iz)v_i$, see
\cite{NV} 2.4. We need to prove an equality 
$c_X(v\otimes zx)=c_X(\sum_i\eps(h_iz)v_i\otimes x)$, or equivalently 
$$\sum_ih_izx\otimes v_i=\sum_i\eps(h_i^{(1)}z)h_i^{(2)}x\otimes v_i$$
which follows from the identity $\eps(h_i^{(1)}z)h_i^{(2)}=h_iz$, see
\cite{NV} Proposition 2.2.1 (v).

To prove that the map $c_X$ is a morphism of $R-$bimodules, we need to
identify $H_s$ and $H_t^{op}$. For this we will use $\eps_s|_{H_t}$ and
$\eps_t|_{H_s}$ (in the case of weak Hopf algebras this coincides with the
usual identification via the antipode).

(2) The map $c_X$ is a map of $H_t-$modules: this reduces to the equality
$$\sum_izh_ix\otimes v_i=\sum_i\eps(\eps_s(z)h_i^{(1)})h_i^{(2)}x\otimes v_i$$
for $z\in H_t, v\in V, x\in X$
which is a consequence of the known identity $zh=\eps(\eps_s(z)h^{(1)})h^{(2)}$
for $h\in H, z\in H_t$.

(3) The map $c_X$ is a map of $H_s-$modules: this reduces to the equality
$$\sum_i\eps(h_i^{(1)}\eps_t(z))h_i^{(2)}x\otimes v_i=\sum_ih_izx\otimes v_i$$
for $z\in H_s, v\in V, x\in X$ which is a consequence of the known identity
$hz=\eps(h^{(1)}\eps_t(z))h^{(2)}$ for $h\in H, z\in H_s$. 
 
The Lemma and the Theorem are proved. \sq

{\bf Remark 12.} There is an elementary description of the category $\cC^*$ in 
terms of the algebra $A$: the category $\cC^*$ is equivalent to the category of
the bimodules over $A$ with $\otimes_A$ as a tensor product (note that this
tensor product is exact since the categories of left and right $A-$modules
are semisimple). Furthermore if the category $\cC$ is braided then this is
the same as category of left modules over $A\otimes A^{op}$ where $A^{op}$ is
the opposite algebra of $A$ (we use the brading to define $A^{op}$ and the
multiplication in the  tensor product of algebras). This shows that the tensor
category of $A-$bimodules depends only on Morita equivalence class of $A$
and sometimes allows one to describe explicitly simple objects of $\cC^*$.

\subsection{Ocneanu cells} It is convenient to say that the category $\cM$
is a right module category over the tensor category $\cC^*_{op}$ (the category
$\cC^*_{op}$ is the same as the category $\cC^*$ but the tensor product is
different: $X\otimes_{op}Y:=Y\otimes X$) and consider the definition of 
$\cC^*$ as associativity constraint $(X\otimes M)\otimes F\to X\otimes (M
\otimes F)$ for $X\in \cC, M\in \cM, F\in \cC^*_{op}$.So we have three 
categories $\cC, \cM, \cC^*_{op}$, four bifunctors 
$\cC\times \cC \to \cC, \cC\times \cM \to \cM, \cM\times \cC^*_{op}\to \cM, 
\cC^*_{op}\times \cC^*_{op}\to \cC^*_{op}$, five associativity constraints
and six hexagon axioms. This situation was axiomatized by A.~Ocneanu and
the corresponding structure constants are known under the name ``Ocneanu
cells'', see e.g. \cite{PZ}. So we consider our formalism above as a 
coordinate free version of Ocneanu cells.

\subsection{Dynamical twists in group algebras} Let $G$ be a finite group
and $A\subset G$ be an abelian subgroup. Then $\Rep(A)$ is a module category
over $\Rep(G)$ and so we have the $R-$fiber functor $F$ where $R$ is a group
algebra of $A$. In \cite{EN} P.~Etingof and D.~Nikshych studied all possible
structures of tensor functor on the functor $F$ (= dynamical twists of the
corresponding weak Hopf algebra up to gauge equivalence) and showed that
they are classified by the isomorphism classes of ``dynamical data'', see
\cite{EN}, Theorem 6.6 for precise statement and examples. From our point
of view they studied module categories $\cM$ over $\Rep(G)$ such that
$K_0(\cM)=K_0(\Rep(A))$ as based modules over $K_0(\Rep(G))$. Using Theorem 2
we see that this is equivalent to looking for pairs $(H,\omega)$ such that 
simple objects of $\Rep^1(H,\omega)$ are numbered by $A^*:=\Hom(A,\BC^*)$ and 
such that $\dim \Hom_{\tilde H}(X\otimes V_\chi, V_\psi)=\dim \Hom_A(X\otimes 
\chi,\psi)$ for any $X\in \Rep(G),\; \chi, \psi \in A^*$ where $V_\chi, V_\psi$
are simple objects in $\Rep^1(H,\omega)$ corresponding to $\chi$ and $\psi$.
Using the Frobenius reciprocity we see that this is equivalent to having
isomorphisms $\Ind_H^G(V_\psi \otimes V_\chi^*)=\Ind_A^G(\psi \chi^{-1})$
of $G-$modules for any $\chi, \psi \in A^*$. But this is precisely the
definition of dynamical data from \cite{EN}. So we see that Theorem 6.6 of 
\cite{EN} is essentially a special case of Theorem 2 above. So we consider 
Theorem 2 as a generalization of \cite{EN}.

\section{Results of B\"ockenhauer, Evans and Kawahigashi}
Let $\cC$ be a rigid monoidal category. In this section we assume in addition
that $\cC$ is {\em braided}, see e.g. \cite{BK}. In this case 
B\"ockenhayer, Evans and Kawahigashi (inspired by A.~Ocneanu) proved a number 
of remarkable results, see \cite{BEK} and references therein. They worked in 
the realm of Operator Algebra Theory. We are
going to translate their results to a categorical language. We will not give
proofs (but note that proofs in \cite{BEK} are calculations in the weak Hopf
algebra attached to a module category over a monoidal category and so should 
work in the general categorical situation), so the skeptical reader may 
consider Claims below as Conjectures (but these Conjectures are Theorems in a 
case when monoidal and module categories under consideration appear in the 
context of subfactor theory; this class is large enough, for example it 
includes all fusion categories of $\widehat{sl}(n)$ thanks to the work of 
A.~Wasserman, see \cite{W}). In any case categorical proofs of statements 
below would be highly desirable.

\subsection{$\alpha-$induction}
Let $\cM$ be a module category over $\cC$ and let $\cC^*$ be the corresponding
dual category, see \ref{dual}. We will assume that all the categories $\cC$,
$\cM$, $\cC^*$ are semisimple, the unit object of $\cC$ is irreducible, and
the category $\cM$ is indecomposable (so the unit object of $\cC^*$ is also
irreducible). Assume that category $\cC$ is braided; let 
$\beta_{V,W} :V\otimes W$ denote braiding and let $\tilde \beta_{W,V}$ denote 
opposite braiding (that is $\beta \tilde \beta =Id$).
We have two tensor functors $\alpha^+, \alpha^-: \cC \to \cC^*$ defined as
follows. For any $X\in \cC$ we have $\alpha^+(X)=\alpha^-(X)=V\otimes ?$ as
functors but module functors structures are different: we set $c^{\alpha^+(X)}_
{Y,M}=\beta_{X,Y}\otimes \id : \alpha^+(X)(Y\otimes M)=X\otimes Y\otimes M\to
Y\otimes X\otimes M=Y\otimes \alpha^+(X)(M)$ and $c^{\alpha^-(X)}_{Y,M}=\tilde
\beta_{X,Y}\otimes \id$. One checks immediately that this defines a module
functors structure using the hexagon axiom. Moreover, $\alpha^+$ and $\alpha^-$
are tensor functors again thanks to the hexagon axiom (here a tensor structure
on $\alpha^\pm$ is given by the associativity constraint: $\alpha^\pm (X\otimes
Y)=X\otimes Y\otimes ?=\alpha^\pm (X)\circ \alpha^\pm (Y)$).

Let $\cC^*_+$ (resp. $\cC^*_-$) denote the additive subcategory of $\cC^*$ 
whoose objects are subquotients ($=$ direct summands) of $\alpha^+(X)$ (resp. 
$\alpha_-(X)$) for all $X\in \cC$. Clearly $\cC^*_+$ and $\cC^*_-$ are
monoidal subcategories of $\cC^*$.
One checks easily that the braiding $\beta_{X,Y}$ defines an isomorphism of
module functors $\alpha^+(X)\circ \alpha^-(Y)\simeq \alpha^-(Y)\circ 
\alpha^+(X)$. 

{\bf Proposition 4.} The braiding above restricts to a well defined functorial
``relative'' brading $\beta^*_{F,G}: F\circ G\to G\circ F$ for all $F\in 
\cC^*_+, G\in \cC^*_-$.

{\bf Proof.} Let $F\in \cC^*$ be a module functor. The module structure on $F$
defines a morphism of functors $F\circ \alpha^\pm (X)\to \alpha^\pm (X)\circ
F$. We will say that $F$ commutes with $\alpha^\pm (X)$ if this morphism is
a morphism of module functors. Clearly if $F$ commutes with $\alpha^\pm (X)$
then the same is true for any subquotient of $F$. One verifies easily that
each functor of the form $\alpha^+(X)$ commutes with any functor of the
form $\alpha^-(Y)$ and vice versa. So the isomorphism $\alpha^+(X)\circ
\alpha^-(Y)\simeq \alpha^-(Y)\circ \alpha^+(X)$ is functorial in the
variables $\alpha^+(X)$ and $\alpha^-(Y)$ (that is commutes with any
endomorphism of these functors) whence we get the Proposition. \sq

In particular consider an additive subcategory $\cC^*_0:=\cC^*_+\cap \cC^*_-
\subset \cC^*$. By the Proposition it has a structure of braided category.

\subsection{Modular invariants} Assume in addition that $\cC$ is a ribbon
category. In this case one defines (see e.g. \cite{BK}) the operators $S$ and 
$T$ acting on the complexified Grothendieck group $K_0(\cC)\otimes \BC$. Let 
$Irr(\cC)=\{ \lambda, \mu, \ldots \}$ be the set indexing simple objects in 
$\cC$. Then $S$ and $T$ can
be considered as matrices with rows and columns indexed by $\lambda, \mu,
\ldots$ by using the basis $[X_\lambda], [X_\mu], \ldots$ of $K_0(\cC)$ 
consisting of classes of simple objects $X_\lambda , X_\mu, \ldots$. Consider 
the matrix $Z_{\lambda, \mu}=\dim \Hom_{\cC^*}(\alpha^+(X_\lambda), 
\alpha^-(X_\mu ))$.

{\bf Claim 1.} The matrix $Z_{\lambda ,\mu}$ commutes with the matrices
$S$ and $T$.

The matrix $Z_{\lambda,\mu}$ evidently has the properties $Z_{\lambda,\mu}\in
\BZ_{\ge 0}$ and $Z_{0,0}=1$ where $0\in Irr(\cC)$ is the index corresponding 
to the trivial object of $\cC$. Such matrices are called {\em modular 
invariants} (under assumption that $S$ is invertible) and were extensively 
studied, see e.g. \cite{Ga} and references therein.

\subsection{Nondegenerate case} Assume in addition that the category $\cC$ is
modular, that is the matrix $S$ is nondegenerate. In this case additional 
results can be established.

{\bf Claim 2.} The category $\cC^*$ is generated by $\cC^*_+$ and $\cC^*_-$,
that is any object of $\cC^*$ is a subquotient of $\alpha^+(X)\circ \alpha^-
(Y)$ for some $X,Y\in \cC$.

{\bf Claim 3.} The number of irreducible objects in the category $\cC^*$ is
$Tr(ZZ^t)=\sum_{\lambda,\mu}Z_{\lambda,\mu}^2$. Moreover the Grothendieck
ring $K_0(\cC^*)$ is isomorphic to the direct sum of matrix algebras of sizes
$Z_{\lambda,\mu}$.

Recall that the characters ($=$ homomorphisms to $\BC$) of $K_0(\cC)$ are 
naturally labelled by $\lambda \in Irr(\cC)$ via the matrix $S$, see \cite{BK}.
The abstract $K_0(\cC)\otimes \BC-$module $K_0(\cM)\otimes \BC$ is a direct 
sum of one dimensional modules, and we refer to the corresponding 
(multi)subset of $Irr(\cC)$ as to the (multi)set of {\em exponents} of $\cM$. 
On the other hand, the multiset containing $\lambda$ with multiplicity 
$Z_{\lambda, \lambda}$ is called the set of exponents of the modular 
invariant $Z$.

{\bf Claim 4.} The number of irreducible objects in the category $\cM$ is
$Tr(Z)=\sum_\lambda Z_{\lambda,\lambda}$. Moreover, the set of exponents of
$\cM$ coincides with the set of exponents of $Z$.

\subsection{The centers} In this section we translate to our language some
results of \cite{BE}.
Recall that there exists a semisimple algebra $A\in \cC$
such that $\cM$ is equivalent to $\Mod_\cC (A)$.

{\bf Definition 15.} The {\em left center} of $A$ is a maximal subobject 
$B^+\subset A$ such that the diagram below commutes:
$$
\xymatrix{B^+\otimes A\ar[rr]^{\beta_{B^+,A}} \ar[dr]^{m} 
&&A\otimes B^+
\ar[dl]_{m}\\ &A&}
$$
The {\em right center} of $A$ is a maximal subobject $B^-\subset A$ such that
the diagram commutes:
$$
\xymatrix{B^-\otimes A\ar[rr]^{\tilde \beta_{B^-,A}} \ar[dr]^{m} 
&&A\otimes B^-
\ar[dl]_{m}\\ &A&}
$$

It is clear that both $B^+$ and $B^-$ are well defined. Note that there is no
reason for the equality $B^+=B^-$ in general (and there are examples when
$B^+$ is not isomorphic to $B^-$ even as an object of $\cC$!). Also it is 
clear that $B^\pm$ are commutative subalgebras of $A$. It is not difficult
to give a definition of $B^\pm$ in terms of the category $\cM$: for 
example, $B^+$ is the universal object $B\in \cC$ endowed with a functorial
morphism $z_M:B\otimes M\to M$ for any $M\in \cM$ such that the diagram
$$\xymatrix{(B\otimes X)\otimes M\ar[rr]^{m_{B,X,M}}
\ar[d]_{\beta_{B,X}\otimes id}&&B\otimes(X\otimes M)\ar[d]^{z_{X\otimes M}}\\
(X\otimes B)\otimes M\ar[rr]^{(id\otimes z_M)\circ m_{X,B,M}}&&X\otimes M}$$
commutes for any $M\in \cM, X\in \cC$. This means that the notion of the
centers is Morita invariant (that is depends only on the class of Morita
equivalence of algebra $A$).

Let $\cM^\pm =\Mod_\cC(B^\pm)$ and let $Z^\pm$ denote the corresponding 
modular invariants. Note that $Z^\pm$ are ``type I modular invariants'', that 
is considered as hermitian forms on $K_0(\cC)$ they are sums of squares of 
linear combinations of characters, see \cite{KO}. Furthermore let 
$\Mod^0_\cC (B^\pm)$ denote the tensor category of representations of the
``vertex algebra'' $B^\pm$ (see \cite{KO}, this category was denoted $Rep^0
(B^\pm)$ there). 

{\bf Claim 5.} We have tensor equivalences $\Mod^0_\cC(B^+)=\cC^*_0=\Mod^0_\cC
(B^-)$. The modular invariants $Z^\pm$ are ``type I parents'' of the modular 
invariant $Z$. In particular the structure of $B^+$ as an object of $\cC$ is 
given by the first (vacuum) row of $Z$ and the structure of $B^-$ as an object
of $\cC$ is given by the first column of $Z$.

We refer the reader to \cite{BE, DMS} for the discussion of type I and type II
modular invariants and the notion of type I parents (due to G.~Moore and 
N.~Seiberg).

\subsection{Problem} We would like to close this section by the following
problem. Let $\lambda \mapsto \bar \lambda$ be the involution of $Irr(\cC)$
induced by the duality: $X_\lambda^*\simeq X_{\bar \lambda}$. It is well
known that the matrix $Z_{\lambda,\mu}=\delta_{\lambda,\bar \mu}$ is a
modular invariant (``charge conjugation'' invariant, see e.g. \cite{DMS}). 

{\bf Problem.} What is a construction of a module category corresponding to
this modular invariant?

Since the charge conjugation modular invariant exists quite generally, one
should expect that there exists a very general construction of this kind.
It is clear from the discussion above that the number of simple objects
in this module category should be equal to the number of selfdual simple
objects in $\cC$. 

\section{Module categories over fusion category of $\widehat{sl}(2)$}
The celebrated result of Capelli-Itzykson-Zuber \cite{CIZ} and Kato \cite{K}
states that $\widehat{sl}(2)-$modular invariants are classified by simply laced
Dynkin diagrams. On the other hand we know that modular invariants can be
constructed from module categories. The aim of this section is to prove
that classification of module categories over fusion categories of $\widehat
{sl}(2)$ is exactly the same: indecomposable module categories are classified
by simply laced Dynkin diagrams.  

\subsection{Monoidal category $\cC_l$}\label{shurik}
Let $l$ be an positive integer. Let $\cC_l$ be the category of 
representations of $\widehat{sl}(2)$ on the level $l$, see e.g. \cite{Kac}.
This category has a natural structure of a monoidal category (fusion 
product), see e.g. \cite{Fi}. Moreover the category $\cC_l$ is modular
category (= braided, balanced, with invertible $S-$matrix), see e.g. 
\cite{BK}. The category $\cC_l$ is semisimple and has $l+1$ simple objects 
denoted by $V_0, V_1, \ldots V_l$ where subscript is highest weight. The 
object $V_0$ is the unit object and the structure of $K_0(\cC_l)$ is 
completely determined by rules
$$[V_1][V_i]=[V_i][V_1]=[V_{i-1}]+[V_{i+1}],\; 1\le i<l;\; [V_1][V_l]=[V_l]
[V_1]=V_{l-1}.$$
In particular the ring $K_0(\cC_l)$ is generated by $[V_1]$. The relations
above imply that $[V_l][V_l]=[V_0]$, that is $V_l$ is an invertible object
(or ``simple current'' in physical language).

It is proved by T.~Kerler, see \cite{Ke, KW}, that there are exactly two 
monoidal categories with the Grothendieck ring isomorphic to $K_0(\cC_l)$, 
first is $\cC_l$ itself and second is a twisted version $\cC_l^{tw}$ of 
$\cC_l$ (the twist comes from the $\BZ/2\BZ$ grading of $K_0(\cC_l)$).

The subcategory of $\cC_l$ additively generated by the objects $V_0$ and $V_l$
is the monoidal subcategory of $\cC_l$ with the Grothendieck ring isomorphic 
to $K_0(\Rep(\BZ/2\BZ))$. The structure of this category is determined by the 
following

{\bf Lemma 7.} (A.~Kirillov, Jr.) The monoidal category generated by $V_0$ and
$V_l$ is equivalent to $\Rep(\BZ/2\BZ)$ for even $l$ and to 
$\Rep(\BZ/2\BZ)^{tw}$ for odd $l$.

{\bf Proof.} The braiding morphism $\beta: V_l\otimes V_l\to V_l\otimes V_l$
is just a number since $\Hom (V_0,V_0)$ is one dimensional. This number
was computed in \cite{KO}, Lemma 7.6 and it equals to $e^{3\pi il/2}$. So this
number equals to $\pm 1$ for even $l$ and $\pm i$ for odd $l$. Now Example
2.5.3 of \cite{Q} implies the Lemma (or better we suggest the reader to deduce
the Lemma from the Hexagon axiom as an useful exercise). \sq

{\bf Remark 13.} It is easy to see that the subcategory generated by $V_0$ and 
$V_l$ in $\cC_l^{tw}$ is always equivalent to $\Rep(\BZ/2\BZ)$ (since for
even $l$ the twist does not touch associativity morphisms in this category
and for odd $l$ the twist untwists the category $\Rep(\BZ/2\BZ)^{tw}$ back
to $\Rep(\BZ/2\BZ)$. 

\subsection{Based modules over $K_0(\cC_l)$} \label{ekdz}
Indecomposable based modules
over $K_0(\cC_l)$ were classified by P.~Di~Francesco and J.-B.~Zuber, see
\cite{DFZ}, and independently by P.~Etingof and M.~Khovanov, see 
\cite{EK}. Such modules are in one to one correspondence with simply laced
Dynkin diagrams (possibly with loops) with Coxeter number $h=l+2$. The 
correspondence is given as follows: let $I$ be the set of vertices of a Dynkin
diagram and consider the matrix $A=(a_{ij})_{i,j\in I}$ given by
$a_{ij}=\mbox{number of edges joining vertices}\; i\; \mbox{and}\; j$ (in 
particular $a_{ii}=0$ for vertices $i$ without loops). Consider a free 
$\BZ-$module $M$ with a basis labelled by $I$. Then letting $[V_1]$ act on $M$
via the matrix $A$ we get a well defined structure of a based module over 
$K_0(\cC_l)$ on $M$.

Recall (see e.g. \cite{EK}) that there are three series of simply laced
Dynkin diagrams with loops: $A_n, D_n$ and $T_n$ (``tadpole'', the only 
diagram with loops) and three exceptionals $E_6, E_7, E_8$. The Coxeter numbers
of these diagrams are respecively $n+1, 2n-2, 2n+1, 12, 18, 30$. We are going 
to prove that for each type except $T_n$ there exists a unique module category
over $\cC_l$ which categorifies the corresponding based module. For type $T_n$
such a module category does not exist.

\subsection{Classification of module categories over $\cC_l$} We will say that
an indecomposable module category $\cM$ over $\cC_l$ is of the type $A_n, D_n,
E_n, T_n$ if the based module $K_0(\cM)$ over $K_0(\cC_l)$ is of the type 
$A_n, D_n, E_n, T_n$ via the correspondence in \ref{ekdz}. In this section we 
are going to prove the following 

{\bf Theorem 6.} For a simply laced Dynkin diagram $X$ of type $A, D, E$ with 
Coxeter number $h$ there exists a unique module category of type $X$ over
$\cC_{h-2}$. The module category of type $T_n$ over $\cC_{2n-1}$ does not
exist.

{\bf Proof.} First we give a construction of all the module categories:

(1) Type $A_n$. This module category is just the ``regular representation'' of
$\cC_l$ that is $\cC_l$ itself considered as a module category over itself.

(2) Type $D_n$. Let $l$ be an even number greater than 2.
By Lemma 7 we have that the subcategory of $\cC_l$ with irreducible objects
$V_0, V_l$ is equivalent to $\Rep(\BZ/2\BZ)$, in particular the object
$A=V_0\oplus V_l$ has a structure of semisimple algebra. One verifies easily
that the corresponding module category $\Mod_{\cC_l}(A)$ has Grothendieck
group isomorphic to based module of type $D_{l/2+2}$ (see e.g. \cite{KO} for
the case of $l=4m$; the case $l=4m+2$ is completely analogous).

(3) Types $E_6, E_8$. The conformal embedding $\widehat{sl}(2)_{10}\subset
\widehat{sp}(4)_1$ determines a structure of a semisimple algebra on the object
$A=V_0\oplus V_6$ of the category $\cC_{10}$, and one verifies easily that 
 the corresponding module category $\Mod_{\cC_{10}}(A)$ has Grothendieck
group isomorphic to the based module of type $E_6$, see \cite{KO}. Similarly,
the conformal embedding $\widehat{sl}(2)_{28}\subset (\widehat{G_2})_1$
defines a structure of a semisimple algebra on the object $A=V_0\oplus V_{10}
\oplus V_{18}\oplus V_{28}$ and the corresponding module category has
Grothendieck group isomorphic to the based module of type $E_8$, see {\em loc.
cit.}

(4) Type $E_7$. In this case $l=16$ and there is no conformal embedding of
$\widehat{sl}(2)_{16}$ in any other affine Lie algebra. Instead there exists
a conformal embedding $\widehat{sl}(2)_{16}\oplus \widehat{sl}(3)_6\subset
(\widehat{E_8})_1$, see \cite{DMS}, Chapter 17. The category of 
representations of $\widehat{sl}(2)_{16}\oplus \widehat{sl}(3)_6$ is equivalent
to the ``tensor product of tensor categories'' $\cC=\cC_{16}\otimes 
\cC(\widehat{sl}(3)_6)$ where $\cC(\widehat{sl}(3)_6)$ is a category of
integrable $\widehat{sl}(3)-$modules on level 6 (we leave to the reader the
definition of tensor product of monoidal categories and the proof of this
statement). The conformal embedding above defines a semisimple algebra $A$
in the category $\cC$; see \cite{DMS} 17.109 for the structure of $A$ as an
object of $\cC$. Next one calculates that the corresponding module category
$\Mod_{\cC}(A)$ has 24 irreducible objects. There are two ways to perform
this calculation: one using explicit fusion rules for $\widehat{sl}(2)_{16}$
and $\widehat{sl}(3)_6$ similarly to the calculations for $D_n, E_6, E_8$, but
this is difficult since the category $\cC$ has $17\cdot 28=476$ irreducible
objects. The second way is to use Claim 4 above (the results of \cite{W}
allow one to apply the results of \cite{BEK} to the category $\cC$) which gives
the desired result from one look at \cite{DMS} 17.109); we prefer this second
way but would like to stress that the usage of Claim 4 can be avoided. 
Now the category
$\cC$ contains the monoidal subcategory $\cC_{16}\otimes \be$ which is
equivalent to $\cC_{16}$. This implies that the category $\cC_{16}$ has a
module category with 24 irreducible objects. Now the results of \ref{ekdz}
shows that any indecomposable module category over $\cC_{16}$ has either
17 (type $A_{17}$) or 10 (type $D_{10}$) or 7 (type $E_7$) irreducible 
objects. Since there are only two decompositions $24=17+7=10+7+7$, we see
immediately that a module category over $\cC_{16}$ of type $E_7$ does exist.
We note that the modular invariants arguments (see \cite{DMS} 17.6) or
explicit calculation show 
that $\cC_{16}-$module category $\Mod_{\cC}(A)$ is a direct sum of the 
category of type $D_{10}$ and two categories of type $E_7$. 

Now we prove that to any simply laced Dynkin diagram with
loops corresponds at most one module category. Let $M$ be the object of
the module category corresponding to the end of the longest leg of the Dynkin 
diagram. It is easy to calculate the structure of $A=\iHom(M,M)$ as an object 
of $\cC_l$. We get the following table from \cite{KO}:

$$\begin{array}{|c|c|c|}\hline \mbox{Diagram}& l=h-2& A\\
\hline A_n&n-1&V_0\\ \hline D_n&2n-4&V_0\oplus V_l\\ \hline T_n&2n-1&V_0\oplus
V_l\\ \hline E_6&10&V_0\oplus V_6\\ \hline E_7&16&V_0\oplus V_8\oplus V_{16}\\
\hline E_8&28&V_0\oplus V_{10}\oplus V_{18}\oplus V_{28}\\ \hline
\end{array}$$

It is immediately clear that $A$ has only one structure of an algebra for type
$A_n$ and no more than one structure of a semisimple algebra for types $D_n, 
T_n$.  

By Lemma 7 the object $A=V_0\oplus V_l$ of $\cC_l$ has a unique 
structure of a semisimple associative algebra for even $l$ and has no such 
structure for odd $l$. In particular, the module category of type $T_n$ 
does not exist and the module category of type $D_n$ does exist and is unique.

To prove uniqueness for types $E_6, E_7$ we will use the following

{\bf Lemma 8.} Let $\cC$ be a rigid monoidal category and $X$ be an 
irreducible 
object of $\cC$. Assume that $\Hom(X\otimes X,X)$ is one dimensional. 
Then $A=\be \oplus X$ has at most one structure of a semisimple
algebra in $\cC$.

{\bf Proof} 
is a word by word repetition of the argument in \cite{KO}, page 24
Type $E_6$. \sq

The Lemma implies immediately that the algebra $A=V_0\oplus V_6$ of type $E_6$
is unique. The algebra of type $E_7$ contains a subalgebra
$A'=V_0\oplus V_{16}$. So we can consider this algebra as an algebra in the 
monoidal category of $A'-$bimodules, $A=A'\oplus X$ where $X$ is $A'-$bimodule
and $X=V_8$ as an object of $\cC_{16}$. One verifies easily along the lines of
\cite{KO} Section 8 that $X$ has four possible structures of an $A'-$bimodule 
and is irreducible as an $A'-$bimodule. For two of these structures 
$X\otimes_{A'}X$ does not contain $A'$ as a direct summand and so these 
bimodules cannot appear in a semisimple algebra $A$. Two other bimodule 
structures are permuted by the automorphism of the algebra $A'$ which is 1 on 
$V_0$ and $(-1)$ on $V_{16}$ so it is enough to consider only one such 
structure. Now we can apply the Lemma above to get the uniqueness of an 
algebra of type $E_7$.

Now consider the case of the algebra $A$ of type $E_8$. In this case again 
$A=A'\oplus X$ where $A'=V_0\oplus V_{28}$ is a subalgebra and $X=V_{10}\oplus
V_{18}$ is an $A'-$bimodule. One shows that there are two possible structures
$X^{(1)}, X^{(2)}$ of an $A'-$bimodule on $X$: $X^{(1)}$ comes from 
$X=\alpha^+(V_{10})$ (where the $\alpha-$induction is taken with respect to 
the category of $A'-$bimodules) and the right $A'-$action of $A'$ on $X^{(2)}$ 
differs from that on $X^{(1)}$ by an automorphism of $A'$ which is 1 on $V_0$ 
and -1 on $V_{28}$ (the left $A'-$actions on $X^{(1)}$ and $X^{(2)}$ are
assumed to be the same). One verifies easily that $X^{(1)}\otimes_{A'}X^{(1)}=
X^{(2)}\otimes_{A'}X^{(2)}$ as $A'-$bimodules and moreover this tensor
product contains $X^{(1)}$ with multiplicity 2 as a direct summand and does
not contain $X^{(2)}$. This implies that in the algebra $A$ one has 
$X=X^{(1)}$ (otherwise the product $X\times X\to X$ is zero which is possible 
iff $X\otimes_{A'}X=A'$) and $A'$ lies in the ``two sided'' center 
(intersection of $B^+$ and $B^-$) of $A$. Now one shows easily that the 
multiplication in $A$ is commutative: the maps $V_{10}\otimes V_{10}\to A$ and
$V_{18}\otimes V_{18}\to A$ commute with braiding by Lemma 7.5 of \cite{KO} 
and the maps $V_{10}\otimes V_{18}\to A$ and $V_{18}\otimes V_{10}\to A$ are 
permuted by the brading thanks to the associativity since $V_{18}$ is the image
of the multiplication of $V_{10}$ and $V_{28}$. Finally, the uniqueness of the
{\em commutative} algebra $A=V_0\oplus V_{10}\oplus V_{18}\oplus V_{28}$ was 
shown in \cite{KO}. Note that using the known structure of the modular 
invariant of type $E_8$ and Claim 5 above one deduces the commutativity of the
algebra of type $E_8$ immediately.

The Theorem is proved. \sq  

{\bf Remark 14.} (i) The Theorem above is not new (but probably our way to 
state it is new). Unfortunately the history is a bit complicated and I don't 
know who should get credit for it. This Theorem was undoubtedly known to 
physicists for some time; I believe that A.~Ocneanu was the first who 
translated this result into mathematical language (in a context of the 
subfactors theory). Unfortunately his results are difficult to understand for 
a nonexpert in Operator Algebras; one of the purposes of this paper is to make
these remarkable results accessible to a student with standard background in 
algebra. At least a big portion of Theorem 6 is contained in \cite{BEK}, in 
particular the subfactor theory construction of the module category of type 
$E_7$, see \cite{BEK} Appendix (unfortunately I don't understand this 
construction; in fact this paper grew up from my attempts to understand it). 
The idea of using the conformal inclusions to construct module categories was 
first translated into mathematical language by F.~Xu, see \cite{X}.  

(ii) Using the Theorem above and Remark 13 one finds
that the indecomposable module categories over $\cC_l^{tw}$ are classified
by the Dynkin diagrams of types $A_n, D_n, T_n, E_6, E_7, E_8$. This gives
some explanation why it is difficult to rule out tadpoles $T_n$ using just
combinatorial methods.  

(iii) The situation with module categories (or modular invariants) of type
$D_n$ has a vast generalization known under the name of simple currents. This
can be summarized as follows: let $\cC$ be a monoidal category such that each
irreducible object of $\cC$ is invertible. Such a category is completely 
determined by the group $G$ of isomorphism classes of irreducible objects
and by the class $\omega \in H^3(G,\BC^*)$, see e.g. \cite{Q, TY}. The 
object $A=\oplus_{g\in G}X_g$ where $X_g$ is a representative of the 
isomorphism class $g$ has a structure of a semisimple algebra if and only if
the class $\omega$ is trivial. In this case the possible structures of 
a semisimple algebra on $A$ are classified by $H^2(G,\BC^*)$ (``discrete 
torsion''), see e.g. \cite{FS}.

(iv) A lot of explicit information on the module categories above is 
available in the literature. For example the structure of the categories
$\cC^*$ (e.g. the Grothendieck ring) is known in all cases thanks to the work 
of A.~Ocneanu. This information is usually presented in the form of
the ``Ocneanu graphs'', see the beautiful pictures e.g. in 
\cite{O, Om, PZ}. Moreover, A.~Ocneanu calculated in all cases categories
$\Fun_{\cC_l}(\cM_1,\cM_2)$ where $\cM_1, \cM_2$ are possibly different
module categories over $\cC_l$, see {\em loc. cit.}
For each irreducible object $F\in \cC^*$ one associates the
``twisted partition function'' which is a matrix $a_{ij}:=\dim \Hom(\be,
\alpha^+(V_i)\otimes F\otimes \alpha^-(V_j))$ where $i,j\in 0, 1,\ldots l$.
The paper \cite{Coq} contains the tables of all twisted partition functions.  

(v) Of course there is an obvious problem to generalize the classification
above to the other simple Lie algebras. One can find in \cite{Ga} a good
account of the related combinatorics. I believe that A.~Ocneanu solved the
corresponding problem for $\hat{sl}(3)$ and $\hat{sl}(4)$, see e.g. his
announcement \cite{Om}. It would be extremely interesting to see the details
of his work (using our methods we are probably able to reprove part of his
results, but I don't know how to construct module categories of types
$\cA_n^*, \cD_n^*$ over $\widehat{sl}(3)$).

(vi) Let $X$ be a Dynkin diagram of type $A, D, E$ and let $\cM$ be the
corresponding module category over $\cC_l$. Let $I=\{ 0,1,\ldots l\}$ be the 
set labeling irreducible objects of $\cC_l$ and let $\fA$ be the set labeling
irreducible objects of $\cM$. For any $a\in \fA$ let $M_a$ be the 
corresponding object. For $i\in I, a,b\in \fA$ consider the vector space
$W_{ab}^i:=\Hom(V_i\otimes M_a,M_b)$ (in the terminology of Ocneanu, this is
the space of essential paths from $a$ to $b$ of length $i$, see \cite{O}). 
Using the canonical
morphisms $V_{i+j}\to V_i\otimes V_j$ (where $V_{i+j}=0$ if $i+j>l$) one
defines a multiplication $W_{ab}^i\otimes W_{bc}^j\to W_{ac}^{i+j}$ which
makes the direct sum $\oplus_{a,b,i}W_{ab}^i$ into an associative algebra.
One recognizes immediately that this algebra is exactly the 
Gelfand-Ponomarev preprojective algebra associated to $X$, see e.g. 
\cite{GP}. So this gives an amazing possibility that the module categories 
over $\cC_l$ are related with the quiver varieties. I don't know if it is 
possible to pursue this relation further. Applying a similar construction to 
the module categories over $\widehat{sl}(3), \widehat{sl}(4)$ etc one gets 
a vast generalization of preprojective algebras and perhaps it would be 
interesting to study these objects.

\end{document}